\documentclass[preprint,12pt]{elsarticle}

\usepackage{amssymb}

\journal{}

\begin{document}

\begin{frontmatter}



\title{Renormalized and entropy solutions for the fractional $p$-Laplacian parabolic equation with $L^1$ data}



\author[label1]{Kaimin Teng}
\ead{tengkaimin2013@163.com}
\address[label1]{Department of Mathematics, Taiyuan University of Technology, Taiyuan 030024, PR China}

\author[label2]{Chao Zhang\corref{cor1}}
\ead{czhangmath@hit.edu.cn}
\address[label2]{Department of Mathematics, Harbin Institute of Technology, Harbin 150001, PR China}

\author[label3]{Shulin Zhou}
\ead{szhou@math.pku.edu.cn}
\address[label3]{LMAM, School of Mathematical Sciences, Peking University, Beijing
100871, PR China}

\cortext[cor1]{Corresponding author.}

\begin{abstract}
In this paper we introduce a natural function class and prove the existence and uniqueness of both nonnegative
renormalized solutions and entropy solutions for the fractional $p$-Laplacian parabolic problem with $L^1$ data. And moreover, we
obtain the equivalence of renormalized solutions and entropy
solutions and establish a comparison result.
\end{abstract}


\end{frontmatter}


\section{Introduction}
\label{}

Suppose that $\Omega$ is a bounded domain in $\mathbb{R}^{N}$
with Lipschitz boundary $\partial\Omega$, $T$ is a positive number.
In this paper we study the following nonlinear parabolic problem
\begin{eqnarray}
\left\{
\begin{array}{ll}
\displaystyle u_t+(-\Delta)_p^s u=f
&\textmd{in} \quad \Omega_T\equiv \Omega\times(0,T),\\[2mm]
u=0 &\textmd{in} \quad (\mathbb R^N\backslash \Omega)\times (0,T),\\[2mm]
u(x,0)=u_0(x) &\textmd{in} \quad \Omega,
\end{array}
\right. \label{main}
\end{eqnarray}
where $0<s<1<p<N$ such that $ps<N$ and $(-\Delta)_p^s$ is the fractional $p$-Laplacian operator which, up to normalization factors, is defined as
\begin{eqnarray*}
(-\Delta)_p^s u(x,t):&=&{\rm P.V.} \int_{\mathbb R^N} \frac{|u(x,t)-u(y,t)|^{p-2}(u(x,t)-u(y,t))}{|x-y|^{N+ps}}\,dy\\[2mm]
&=&\lim_{\varepsilon\downarrow 0}\int_{\mathbb R^N\backslash B_\varepsilon(x)} \frac{|u(x,t)-u(y,t)|^{p-2}(u(x,t)-u(y,t))}{|x-y|^{N+ps}}\,dy,
\end{eqnarray*}
where $(x, t)\in \mathbb R^N\times \mathbb R^+$, P.V. is a commonly used abbreviation in the principal value sense. Moreover, we assume that $f$ and $u_0$ are nonnegative satisfying
\begin{equation}
\label{assum}
f\in L^1(\Omega_T) \quad \textmd{and}\quad u_0\in L^1(\Omega).
\end{equation}

There have been a large number of research activities on the study of
well-posedness of $p$-Laplacian type equations and general Leray-Lions problems with $L^1$ and measure data.
Under these assumptions, the existence of a distributional solution, so-called SOLA (Solutions Obtained as Limit of
Approximations), was proved in \cite{BDGL, BG, DA}, but due to the lack of regularity of the solution, the distributional formulation is not
strong enough to provide uniqueness. To overcome this difficulty, it is reasonable to work with renormalized
solutions and entropy solutions, which need less regularity than weak solutions. The notion of renormalized solutions was
first introduced by DiPerna and Lions \cite{DL} for the study of
Boltzmann equation. It was then adapted to the study of some
nonlinear elliptic or parabolic problems and evolution problems in
fluid mechanics. We refer to \cite{BM,BMR,BPR,BR,BGDM,DMOP,LI2,ZZ2,ZZ3} for
details. At the same time the notion of entropy solutions has been
proposed by B\'{e}nilan et al. in \cite{B} for the nonlinear
elliptic problems. This framework was extended to related problems \cite{AB,BC,BGO,LP,P,PR,SU}.

The fractional Laplacian operators and non-local operators have attracted increasing attention over the last years. This type of operators arises in a natural way in many different applications such as continuum mechanics, phase transition phenomena, population dynamics, image process, game theory and L\'{e}vy processes, see for example \cite{A, C1, CS, CV, MK}. For this reason it is particularly important to study situations when such non-local operators
are involved in equations featuring singular or irregular data. This leads to study
non-local equations having $L^1$ or measure data. As far as the non-local $p$-Laplacian operator $(-\Delta)_p^s$ is concerned, the linear elliptic case $p=2$ has been studied in \cite{AAB, KPU, LPPS}. In particular, the existence and uniqueness of renormalized solutions for the problems of the kind
$$
\beta(u)+(-\Delta)^{s}u\ni f \quad\textmd{in } \mathbb R^N
$$
was proved by Alibaud, Andreianov and Bendahmane in \cite{AAB}, where $f\in L^1(\mathbb R^N)$ and $\beta$ is a maximal monotone graph in $\mathbb R$. Using a duality argument, in the sense of Stampacchia, Kenneth, Petitta and Ulusoy in \cite{KPU} proved the existence and uniqueness of solutions to non-local problems like $
(-\Delta)^{s}u=\mu$ in  $\mathbb R^N$ with $\mu$ being a bounded Radon measure whose support is compactly contained in $\mathbb R^N$. In \cite{KMS}, Kuusi, Mingione and Sire discussed the elliptic non-local case $p\not=2$ with measure data and developed an existence of SOLA, regularity and Wolf potential theory.  In addition, Abdellaoui et al in \cite{AAB2} investigated the fractional elliptic $p$-Laplacian equations with weight and general datum and showed that there exists a unique entropy positive solution. On the other hand, Abdellaoui et al in \cite{AABP} established the results on the existence of a weak solution obtained as limit of approximations (SOLA) and the existence of nonnegative entropy solutions for the fractional $p$-Laplacian equations.

In this paper, we focus our attention on the well-posedness of renormalized solutions and the uniqueness of entropy solutions for the fractional $p$-Laplacian parabolic problem (\ref{main}). Our results cover the case of linear parabolic non-local
equations and are also new in such cases for the study of renormalized
solutions. We construct an approximate solution sequence and
establish some $a\ priori$ estimates. Then we draw a subsequence to
obtain a limit function, and prove this function is a renormalized
solution. Based on the convergence results of approximate solutions, we obtain that the renormalized solution of
problem (\ref{main}) is also an entropy solution, which leads to an
inequality in the entropy formulation. By choosing suitable test
functions, we prove the uniqueness of renormalized solutions and
entropy solutions, and thus the equivalence of renormalized
solutions and entropy solutions. Here we would like to mention
that the definition of renormalized solutions is influenced by \cite{AAB}. The main point is to circumvent the use of chain rules, which is
not available in the non-local framework.

For the convenience of the readers, we recall some definitions and
basic properties of the fractional Sobolev spaces, in which main
results can be found in \cite{DGP, LL,MRT,T,XZR} and the references therein.

Let $s\in (0,1)$ and $p>1$. The fractional Sobolev space
$$
W^{s,p}(\mathbb R^N)\equiv \left\{u\in L^p(\mathbb R^N): \int_{\mathbb R^N}\!\!\int_{\mathbb R^N} \frac{|u(x)-u(y)|^p}{|x-y|^{N+ps}}\,dxdy<+\infty\right\}
$$
is a Banach space endowed with the norm
$$
\|u\|_{W^{s,p}(\mathbb R^N)}=\|u\|_{L^p(\mathbb R^N)}+\left(\int_{\mathbb R^N}\!\!\int_{\mathbb R^N} \frac{|u(x)-u(y)|^p}{|x-y|^{N+ps}}\,dxdy\right)^{\frac 1p}.
$$

Denote $D_\Omega=(\mathbb R^N\times \mathbb R^N) \backslash (\mathcal C\Omega\times \mathcal C\Omega)$, where $\mathcal C\Omega=\mathbb R^N\backslash \Omega$.
For every function   $u\in C_0^\infty(\Omega)$ we  define $u=0$ in $\mathcal C\Omega$ and then have  $u\in C_0^\infty(\mathbb R^N)\subset W^{s,p}(\mathbb R^N)$.

Now we define $X_0^{s,p}(\Omega)$ to be the closure of $C_0^\infty(\Omega)$ in $W^{s,p}(\mathbb R^N)$.
For every function  $u\in X_0^{s,p}(\Omega)$, it is clear that
$$
 u=0 \quad \textmd{a.e. in }  \mathcal C\Omega
$$
and
\begin{eqnarray*}
\int_{\mathbb R^N}\!\!\int_{\mathbb R^N} \frac{|u(x)-u(y)|^p}{|x-y|^{N+ps}}\,dxdy
&=&\int_{\Omega}\int_{\Omega} \frac{|u(x)-u(y)|^p}{|x-y|^{N+ps}}\,dxdy\\
&&+2\int_{\Omega}|u(x)|^p\int_{\mathcal C\Omega}  \frac{1}{|x-y|^{N+ps}}\,dydx.
\end{eqnarray*}

Recall  Lemma 6.1 in \cite{DGP}, we have
$$
\int_{\mathcal C\Omega}  \frac{1}{|x-y|^{N+ps}}\,dy\ge c|\Omega|^{-\frac{sp}{N}},
$$
where $c=c(N,p,s)>0$ and then obtain the Poincar\'{e} inequality
\begin{equation}\label{poincare}
\int_{\Omega}|u(x)|^p\,dx\le C\int_{D_\Omega} |u(x)-u(y)|^p\,d\nu, \quad\forall p\ge 1,
\end{equation}
where
$$
d\nu=\frac{dxdy}{|x-y|^{N+ps}}.
$$
Therefore, there exists a positive constant $C=C(N,p,s,\Omega)$ such that for any $u\in X_0^{s,p}(\Omega)$,
$$
\int_{D_\Omega} |u(x)-u(y)|^p\,d\nu \le \|u\|^p_{W^{s,p}(\mathbb R^N)}\le C\int_{D_\Omega} |u(x)-u(y)|^p\,d\nu.
$$

Thus we can endow $X_0^{s,p}(\Omega)$ with the equivalent norm
$$
\|u\|_{X_0^{s,p}(\Omega)}=\left(\int_{D_\Omega}|u(x)-u(y)|^p\,d\nu\right)^{\frac 1p}.
$$
Note that $X_0^{s,p}(\Omega)$ is a uniformly convex Banach space, and hence $X_0^{s,p}(\Omega)$ is a reflexive Banach space.

For $w\in W^{s,p}(\mathbb R^N)$, we define the fractional $p$-Laplacian as
$$
(-\Delta)_p^s w(x)={\rm P.V.} \int_{\mathbb R^N}\frac{|w(x)-w(y)|^{p-2}(w(x)-w(y))}{|x-y|^{N+ps}}\,dy.
$$
It is clear that for all $w, v\in W^{s,p}(\mathbb R^N)$, we have
$$
\langle (-\Delta)_p^s w, v\rangle=\frac 12 \int_{\mathbb R^N\times \mathbb R^N}|w(x)-w(y)|^{p-2}(w(x)-w(y))(v(x)-v(y))\,d\nu.
$$
Now, if $w, v\in X_0^{s,p}(\Omega)$, then
$$
\langle (-\Delta)_p^s w, v\rangle=\frac 12 \int_{D_\Omega} |w(x)-w(y)|^{p-2}(w(x)-w(y))(v(x)-v(y))\,d\nu.
$$

It is easy to check that $(-\Delta)_p^s: X_0^{s,p}(\Omega)\to X_0^{s,p}(\Omega)^*$, where $X_0^{s,p}(\Omega)^*$ denotes the dual space of $X_0^{s,p}(\Omega)$.
Let us define now the corresponding parabolic spaces. As in the local case, the space $L^p(0,T;X_0^{s,p}(\Omega))$ is defined as the set of function $u$ such that $u\in L^p(\Omega_T)$ with $\|u\|_{L^p(0,T;X_0^{s,p}(\Omega))}<\infty$, where
$$
\|u\|_{L^p(0,T;X_0^{s,p}(\Omega))}=\left(\int_0^T\!\!\int_{D_\Omega}|u(x,t)-u(y,t)|^p\,d\nu dt\right)^{\frac 1p}.
$$
$L^p(0,T;X_0^{s,p}(\Omega))$ is a Banach space whose dual space is $L^{p'}(0,T;X_0^{s,p}(\Omega)^*)$.

\medskip

Let $T_k$ denote the truncation function at height $k\ge 0$:
\begin{displaymath}
T_k(r)=\min\{k, \max\{r,-k\}\}=\left\{ \begin{array}{ll} k &
\textrm{if $r\ge k$},\\[2mm]
r & \textrm{if $|r|<k$},\\[2mm]
-k & \textrm{if $r\le -k$},
\end{array} \right.
\end{displaymath}
and its primitive $\Theta_k: \mathbb R\to \mathbb R^+$ by
\begin{displaymath}
\Theta_k(r)=\int_0^r T_k(s)\,ds=\left\{ \begin{array}{ll}
\frac{r^2}{2} &
\textrm{if $|r|\le k$},\\[2mm]
k|r|-\frac{k^2}{2} & \textrm{if $|r|\ge k$}.
\end{array} \right.
\end{displaymath}
It is obvious that $\Theta_k(r)\ge 0$ and $\Theta_k(r)\le k|r|$.

We denote $u\in \mathcal {T}_0^{s,p}(\Omega_T)$ if $u: \mathbb R^N \times (0,T] \to
\mathbb R$ is measurable and $T_k(u)\in
L^{p}\big(0,T;X_0^{s,p}(\Omega)\big)$ for every $k>0$. It is obvious that $u=0$ a.e. in $\mathcal C\Omega$.
For simplicity and for any measurable function $u$, we write
$$
U(x,y,t)=u(x,t)-u(y,t).
$$

Next we give the following
definitions of renormalized solutions and entropy solutions for
problem (\ref{main}).

\bigskip \noindent{\bf Definition 1.1.}
 \label{def2}
A function $u\in \mathcal {T}_0^{s,p}(\Omega_T)\cap
C([0,T];L^1(\Omega))$ is a renormalized solution to problem
{\rm(\ref{main})} if the following conditions are satisfied:

(i) $$\displaystyle \lim\limits_{h\to \infty}\int\!\!\int\!\!\int_{\{(x,y,t):(u(x,t), u(y,t))\in R_h\}}|U(x,y,t)|^{p-1}\,d\nu dt=0,$$ where $$R_h=\Big\{(u,v)\in \mathbb R^{2}: h+1\le \max\{|u|, |v|\} \textmd{ and } (\min\{|u|, |v|\}\le h \textmd{ or } uv<0)\Big\}.$$

(ii) For every function $\varphi\in C^1(\bar\Omega_T)$ with $\varphi=0$ in $\mathcal C\Omega\times (0,T)$ and $\varphi(\cdot,T)=0$ in $\Omega$, and $S\in W^{1,\infty}(\mathbb R)$ which is
piecewise $C^1$ satisfying that $S'$ has a compact support,
\begin{eqnarray}
\label{1-5} &&-\int_\Omega
S(u_0)\varphi(x,0)\,dx-\int_0^T\!\!\int_\Omega S(u)\frac{\partial
\varphi}{\partial t}\,dxdt \nonumber\\
&&\quad+\frac 12\int_0^T\!\!\int_{D_\Omega} |U(x,y,t)|^{p-2}U(x,y,t)[(S'(u)\varphi)(x,t)-(S'(u)\varphi)(y,t)]\,d\nu dt \nonumber\\
&&=\int_0^T\!\!\int_\Omega fS'(u)\varphi\,dxdt
\end{eqnarray}
holds.

\medskip
\noindent{\bf Remark 1.1.} It is not difficult to see that the symmetrization of the difference $(S'(u)\varphi)(x,t)-(S'(u)\varphi)(y,t)$ can yield the following equality:
\begin{eqnarray*}
&&\frac 12\int_0^T\!\!\int_{D_\Omega} |U(x,y,t)|^{p-2}U(x,y,t)[(S'(u)\varphi)(x,t)-(S'(u)\varphi)(y,t)]\,d\nu dt \nonumber\\
&&=\frac 12\int_0^T\!\!\int_{D_\Omega} |U(x,y,t)|^{p-2}U(x,y,t)(S'(u)(x,t)-S'(u)(y,t)) \nonumber\\
&&\qquad\qquad\cdot \frac{\varphi(x,t)+\varphi(y,t)}{2}\,d\nu dt \nonumber\\
&&\quad+\frac 12\int_0^T\!\!\int_{D_\Omega} |U(x,y,t)|^{p-2}U(x,y,t)(\varphi(x,t)-\varphi(y,t)) \nonumber\\
&&\qquad\qquad\cdot\frac{S'(u)(x,t)+S'(u)(y,t)}{2}\,d\nu dt.
\end{eqnarray*}

\medskip \noindent{\bf Definition 1.2.}
 \label{def1}
A function $u\in \mathcal {T}_0^{s,p}(\Omega_T)\cap
C([0,T];L^1(\Omega))$ is an entropy solution to problem
{\rm(\ref{main})} if
\begin{eqnarray}\label{1-4}
&&\int_\Omega \Theta_k(u-\phi)(T)\,dx-\int_\Omega
\Theta_k(u_0-\phi(0))\,dx+\int_0^T \langle \phi_t, T_k(u-\phi)
\rangle\,dt \nonumber\\
&&\quad+\frac 12\int_0^T\!\!\int_{D_\Omega} |U(x,y,t)|^{p-2}U(x,y,t) \nonumber\\[2mm]
&&\qquad\qquad \cdot[T_k(u(x,t)-\phi(x,t))-T_k(u(y,t)-\phi(y,t))]\,d\nu dt\nonumber\\[2mm]
&&\le \int_0^T\!\!\int_\Omega f T_k(u-\phi)\,dxdt,
\end{eqnarray}
for all $k>0$ and $\phi\in C^1(\bar \Omega_T)$ with $\phi=0$ in $\mathcal C\Omega\times (0,T)$.

\bigskip
Now we state our main results. The first two theorems are about the existence and
uniqueness of nonnegative renormalized and entropy solutions. The third one is about the comparison principle.

\medskip
\noindent{\bf Theorem 1.1.} {\it \label{thm1} Assume that condition
{\rm (\ref{assum})} holds. Then there exists a unique renormalized
solution for problem {\rm (\ref{main})}.}

\medskip
\noindent{\bf Theorem 1.2.} {\it \label{thm2}
 Assume that condition {\rm (\ref{assum})} holds. Then the renormalized
solution $u$ obtained in Theorem  {\rm 1.1} is also an entropy solution for
problem {\rm (\ref{main})}. And the entropy solution is unique.}

\medskip
\noindent{\bf Remark 1.2.} The renormalized  solution for problem
{\rm(\ref{main})} is equivalent to the entropy solution for problem
{\rm(\ref{main})}.

\medskip
\noindent{\bf Theorem 1.3.} {\it \label{thm3}
Let $u_0, v_0\in L^1(\Omega)$, $f, g\in L^1(\Omega_T)$
such that $u_0\le v_0$ and $f\le g$. If $u$ is the entropy solution (renormalized solution)
of problem {\rm (\ref{main})} and $v$ is  the entropy solution (renormalized solution) of
problem {\rm (\ref{main})} with  $u_0, f$ being replaced by $v_0,g$,
then $u\le v$ a.e. in $\Omega_T$.}

\bigskip

The rest of this paper is organized as follows. In Section $2$, we
prove the existence and uniqueness of weak solutions to problem (\ref{main}). We will prove  the
main results  in Section $3$. In the following sections $C$ will
represent a generic constant that may change from line to line even
if in the same inequality.

\section{Weak solutions}
\label{} \setcounter{equation}{0}

In this section we will give a reasonable definition for weak solutions and prove
the existence and uniqueness of weak solutions to problem (\ref{main}).

\bigskip

\noindent{\bf Lemma 2.1.} {\it Assume that $u_0\in L^2(\Omega)$ and
$f\in L^{p'}(0,T;X_0^{s,p}(\Omega)^*)$. Then the following problem
\begin{eqnarray*}
\left\{
\begin{array}{ll}
\displaystyle u_t+(-\Delta)_p^s u=f
&\textmd{in} \quad \Omega_T,\\[2mm]
u=0 &\textmd{in} \quad \mathcal C\Omega \times (0,T),\\[2mm]
u(x,0)=u_0(x) &\textmd{in} \quad \Omega
\end{array}
\right.
\end{eqnarray*}
admits a unique weak solution $u\in
L^{p}\big(0,T;X_0^{s,p}(\Omega)\big)\cap
C([0,T];L^2(\Omega))$ with $u_t\in L^{p'}(0,T;X_0^{s,p}(\Omega)^*)$ such
that for any $\varphi\in C_0^\infty(\Omega_T)$,
\begin{eqnarray*}
&&\int_0^T\langle u_t, \varphi\rangle\,dt+\frac 12 \int_0^T\!\!\int_{D_\Omega}|U(x,y,t)|^{p-2}U(x,y,t)(\varphi(x,t)-\varphi(y,t))\,d\nu dt\\
&&=\int_0^T\!\!\int_\Omega f \varphi\,dxdt
\end{eqnarray*}
holds.}

\medskip\noindent{\bf Proof.} Since the fractional $p$-Laplacian operator $(-\Delta)_p^s$ is monotone, the existence of weak solutions can be proved by employing the difference and variation methods. We give a sketched proof.

Let $n$ be a positive integer. Denote $h=T/n$. We first consider the
following time-discrete problem
\begin{eqnarray}\label{2-1}
\left\{
\begin{array}{ll}
\displaystyle \frac{u_k-u_{k-1}}{h}+(-\Delta)_p^s u_k=[f]_h((k-1)h)),\\[3mm]
u_k|_{\mathcal C\Omega}=0,\quad k=1,2,\dots,n,
\end{array}
\right.
\end{eqnarray}
where $[f]_h$ denotes the Steklov average of $f$ defined by
$$
[f]_h(x,t)=\frac 1h\int_t^{t+h}f(x,\tau)\,d\tau.
$$

For $k=1$, we introduce the variational problem
$$
\min\{J(u)|u\in W\},
$$
where
$$
W=\left\{u\in X_0^{s,p}(\Omega)\cap L^2(\Omega)\right\}
$$
and functional $J$ is
$$
J(u)=\frac{1}{2h}\int_\Omega (u-u_0)^2\,dx+\frac{1}{2p}\int_{D_\Omega}
|u(x)-u(y)|^p\,d\nu-\int_\Omega [f]_h(0) u\,dx.
$$
By using the classical Direct Methods of the
Calculus of Variations in fractional Sobolev spaces, we can prove that $J(u)$ is lower bounded and coercive on $W$. On the other
hand, $J(u)$ is  weakly lower semicontinuous on $W$. Therefore,
there exists  a function $u_1\in W$ such that
$$
J(u_1)=\inf\limits_{u\in W}J(u).
$$
Thus the function $u_1$ is a weak solution of the corresponding
Euler-Lagrange equation of $J(u)$, which is (\ref{2-1}) in the case
$k=1$. And it is unique since $J(u)$ is strictly convex.

Following the same procedures, we find weak solutions $u_k$ of
(\ref{2-1}) for $k=2,\dots,n$. It follows that, for every
$\varphi\in W$,
\begin{equation}
\label{2-3} \int_\Omega \frac{u_k-u_{k-1}}{h}\varphi\,dx+\int_\Omega
(-\Delta)_p^s u_k \varphi\,dx=\int_\Omega
[f]_h((k-1)h)\varphi\,dx.
\end{equation}

For every $h=T/n$, we define the approximate solutions
\begin{eqnarray*}
u_h(x,t)=\left\{
\begin{array}{ll}
u_{0}(x), & t=0, \\
u_{1}(x), & 0<t\leq h, \\
\cdots \cdots , & \cdots \cdots , \\
u_{j}(x), & (j-1)h<t\leq jh, \\
\cdots \cdots , & \cdots \cdots , \\
u_{n}(x), & (n-1)h<t\leq nh=T.
\end{array}
\right.
\end{eqnarray*}
Taking $\varphi=u_k$ in (\ref{2-3}), we can obtain an $a\ priori$
estimate
\begin{eqnarray*}
\int_\Omega u^2_h(x,t)\,dx+\frac 12\int_0^T\!\!\int_{D_\Omega} |u_h(x,t)-u_h(y,t)|^{p}\,d\nu dt\le C,
\end{eqnarray*}
which implies that
$$
\|u_h\|_{L^\infty(0,T;L^2(\Omega))}+\|u_h\|_{L^p(0,T;X_0^{s,p}(\Omega))}\le C.
$$
Thus we may choose a subsequence (we also denote it by the original
sequence for simplicity) such that
\begin{eqnarray*}
\begin{array}{ll}
u_h\rightharpoonup u,\quad \textmd{weakly-* in}\quad L^{\infty}(0,T;
L^2(\Omega)), \\[2mm]
u_h\rightharpoonup u,\quad \textmd{weakly in}\quad L^p(0,T;X_0^{s,p}(\Omega)).
\end{array}
\end{eqnarray*}

Recalling the fact that $u\in
L^{p}(0,T;X_0^{s,p}(\Omega))\cap L^\infty(0,T;L^2(\Omega))$ and
$u_t\in L^{p'}(0,T;X_0^{s,p}(\Omega)^*)$, we
conclude that $u$ belongs to $C([0,T];L^2(\Omega))$. Therefore, we
obtain the existence of weak solutions.

For uniqueness, suppose there exist two weak solutions $u$ and $v$
of problem (\ref{main}). Then $w=u-v$ satisfies the following
problem
\begin{eqnarray*}
\left\{
\begin{array}{ll}
\displaystyle w_t+\big[(-\Delta)_p^s u-(-\Delta)_p^s v\big]=0 &\textmd{\rm in} \quad \Omega_T,\\[2mm]
w=0 &\textmd{\rm in} \quad \mathcal C\Omega \times (0,T),\\[2mm]
w(x,0)=0 &\textmd{\rm in} \quad \Omega.
\end{array}
\right.
\end{eqnarray*}
Choosing $w$ as a test function in the above problem, we have, for
almost every $t\in (0,T)$,
\begin{eqnarray*}
\int_\Omega w^2(t)\,dx+\int_0^t\!\!\int_{D_\Omega} [|U(x,y,\tau)|^{p-2}U(x,y,\tau)-|V(x,y,\tau)|^{p-2}V(x,y,\tau)]\\
\cdot (U(x,y,\tau)-V(x,y,\tau))\,d\nu d\tau=0,
\end{eqnarray*}
where $V(x, y, \tau)=v(x,\tau)-v(y,\tau)$. Since the two terms on the left-hand side are nonnegative, then we have
$u=v$ a.e. in $\Omega_T$. This finishes the proof.  $\quad\Box$

\section{The proof of main results}
\label{}

\setcounter{equation}{0}

Now we are ready to prove the main results. Some of the reasoning is
based on the ideas developed in \cite{AABP,PO,PR,ZZ2}. First we  prove the existence and uniqueness
of renormalized solutions for problem (\ref{main}).

\medskip

\noindent{\bf Proof of Theorem 1.1.} (1) Existence of renormalized solutions.

We first introduce the approximate problems. Define $f_n=T_n(f)$ and $u_{0n}=T_n(u_0)$, then we know that $f_n, u_{0n}$ are nonnegative, $(f_n,u_{0n})\in L^\infty(\Omega_T)\times L^\infty(\Omega)$ and $(f_n,u_{0n})\nearrow (f,u_0)$  strongly in $L^1(\Omega_T)\times L^1(\Omega)$ such that
\begin{equation}
\label{assume} \|f_n\|_{L^1(\Omega_T)}\le \|f\|_{L^1(\Omega_T)}, \quad
\|u_{0n}\|_{L^1(\Omega)}\le \|u_0\|_{L^1(\Omega)}.
\end{equation}
Then we consider the approximate problem of (\ref{main})
\begin{eqnarray}
\label{appro} \left\{
\begin{array}{ll} \displaystyle
(u_n)_t+(-\Delta)_p^s u_n=f_n &\textmd{in} \quad \Omega_T,\\[2mm]
u_n=0 &\textmd{in} \quad \mathcal C\Omega\times (0,T),\\[2mm]
u_n(x,0)=u_{0n} &\textmd{in} \quad \Omega.
\end{array}
\right.
\end{eqnarray}
By Lemma 2.1 and comparison principle, we can find a unique nonnegative weak solution $u_n\in
L^{p}\big(0,T; X_0^{s,p}(\Omega)\big)$ for problem (\ref{appro}). Our aim is to prove
that a subsequence of these approximate solutions $\{u_n\}$
converges increasingly to a measurable function $u$, which is a renormalized
solution of problem (\ref{main}). We will divide the proof into
several steps. We present a self-contained proof for the sake of clarity and
readability.

\medskip

\textbf{Step 1.} Prove the convergence of $\{u_n\}$ in
$C([0,T];L^1(\Omega))$ and  find its subsequence which is almost
everywhere convergent in $\Omega_T$.

Let $m$ and $n$ be two integers, then from (\ref{appro}) we can
write the weak form as
\begin{eqnarray*}
&&\int_0^T \langle
(u_n-u_m)_t,\phi\rangle\,dt+\int_0^T\langle(-\Delta)_p^s u_n-(-\Delta)_p^s u_m, \phi\rangle\,dt \\
&&=\int_0^T\!\!\int_\Omega(f_n-f_m)\phi\,dxdt,
\end{eqnarray*}
for all $\phi\in L^{p}(0,T;X_0^{s,p}(\Omega))\cap L^\infty(\Omega_T)$. Choosing
$\phi=T_1(u_n-u_m)\chi_{(0,t)}$ with $t\le T$, we have
\begin{eqnarray*}
&&\int_0^t \langle
(u_n-u_m)_t, T_1(u_n-u_m)\rangle\,d\tau\\
&&\quad+\int_0^t\langle(-\Delta)_p^s u_n-(-\Delta)_p^s u_m, T_1(u_n-u_m)\rangle\,d\tau \\
&&=\int_0^T\!\!\int_\Omega(f_n-f_m)T_1(u_n-u_m)\chi_{(0,t)}\,dxdt.
\end{eqnarray*}
Observe that
\begin{eqnarray*}
&&\int_0^t\langle(-\Delta)_p^s u_n-(-\Delta)_p^s u_m, T_1(u_n-u_m)\rangle\,d\tau\\
&&=\frac12\int_0^t\!\!\int_{D_\Omega}(|U_n(x,y,\tau)|^{p-2}U_n(x,y,\tau)-|U_m(x,y,\tau)|^{p-2}U_m(x,y,\tau))\\[2mm]
&&\qquad\quad \cdot [T_1(u_n(x,\tau)-u_m(x,\tau))-T_1(u_n(y,\tau)-u_m(y,\tau))]\,d\nu d\tau.
\end{eqnarray*}
Since
\begin{eqnarray*}
&&T_1(u_n(x,\tau)-u_m(x,\tau))-T_1(u_n(y,\tau)-u_m(y,\tau))\\[2mm]
&&=T'_1(\xi_{nm})(U_n(x,y,\tau)-U_m(x,y,\tau))
\end{eqnarray*}
due to the mean value theorem, where $T'_1\ge 0$, we know that
$$
\int_0^t\langle(-\Delta)_p^s u_n-(-\Delta)_p^s u_m, T_1(u_n-u_m)\rangle\,d\tau\ge 0.
$$
Then we get
\begin{eqnarray*}
\int_\Omega \Theta_1(u_n-u_m)(t)\,dx&\le& \int_\Omega
\Theta_1(u_{0n}-u_{0m})\,dx+\|f_n-f_m\|_{L^1(\Omega_T)}\\
&\le& \|u_{0n}-u_{0m}\|_{L^1(\Omega)}+\|f_n-f_m\|_{L^1(\Omega_T)}:=a_{n,m}.
\end{eqnarray*}

Therefore, we conclude that
\begin{eqnarray*}
&&\int_{\{|u_n-u_m|<1\}}\frac{|u_n-u_m|^2(t)}{2}\,dx+\int_{\{|u_n-u_m|\ge
1\}}\frac{|u_n-u_m|(t)}{2}\,dx\\
&&\le \int_\Omega [\Theta_1(u_n-u_m)](t)\,dx\le a_{n,m}.
\end{eqnarray*}
It follows that
\begin{eqnarray*}
\int_\Omega |u_n-u_m|(t)\,dx&=&\int_{\{|u_n-u_m|<1\}}
|u_n-u_m|(t)\,dx\\
&&+\int_{\{|u_n-u_m|\ge 1\}} |u_n-u_m|(t)\,dx\\
&\le& \left(\int_{\{|u_n-u_m|<1\}}
|u_n-u_m|^2(t)\,dx\right)^{\frac 12}{\rm meas}(\Omega)^{\frac 12}+2a_{n,m}\\
&\le& (2{\rm meas}(\Omega))^{\frac 12}a_{n,m}^{\frac 12}+2a_{n,m}.
\end{eqnarray*}
Since $\{f_n\}$ and $\{u_{0n}\}$ are convergent in $L^1$, we have
$a_{n,m}\to 0$ for $n,m\to +\infty$. Thus $\{u_n\}$ is a Cauchy
sequence in $C([0,T];L^1(\Omega))$ and $u_n$ converges to $u$ in
$C([0,T];L^1(\Omega))$. Then we find an a.e. convergent  subsequence
(still denoted by $\{u_n\}$) in $\Omega_T$ such that
\begin{equation}
\label{3-4} u_n\nearrow u \quad \mbox{a.e. in } \Omega_T.
\end{equation}

\textbf{Step 2.} Prove $T_k(u_n)$ strongly converges to $T_k(u)$ in $L^p(0,T;X_0^{s,p}(\Omega))$, for every $k>0$.

Choosing $T_k(u_n)$ as a test function in (\ref{appro}), we have
\begin{eqnarray*}
 &&\int_\Omega \Theta_k(u_n)(T)\,dx-\int_\Omega
\Theta_k(u_{0n})\,dx
\\
&&\quad +\int_0^T\langle(-\Delta)_p^s u_n, T_k(u_n)\rangle\,dt
=\int_0^T\!\!\int_\Omega f_n T_k(u_n)\,dxdt.
\end{eqnarray*}
It follows from the definition of $\Theta_k(r)$, $0\le T'_k\le 1$ and (\ref{assume})
that
\begin{eqnarray*}
&&\frac 12\int_0^T\!\!\int_{D_\Omega} |T_k(u_n(x,t))-T_k(u_n(y,t))|^p\,d\nu dt \\
&&\le \frac 12 \int_0^T\!\!\int_{D_\Omega} |U_n(x,y,t)|^{p-2}U_n(x,y,t)(T_k(u_n(x,t)-T_k(u_n(y,t)))\,d\nu dt \\
&&\le k\big(\|f_n\|_{L^1(\Omega_T)}+\|u_{0n}\|_{L^1(\Omega)}\big)
 \\
&&\le k\big(\|f\|_{L^1(\Omega_T)}+\|u_{0}\|_{L^1(\Omega)}\big).
\end{eqnarray*}
Then, up to a subsequence, we deduce that
\begin{eqnarray*}
T_k(u_n)\rightharpoonup T_k(u) \quad\textmd{weakly in } L^p(0,T; X_0^{s,p}(\Omega)).
\end{eqnarray*}

In order to deal with the time derivative of truncations, we will
use the regularization method of Landes \cite{L} and use the
sequence $(T_k(u))_\mu$ as approximation of $T_k(u)$. For $\mu>0$,
we define the regularization in time of the function $T_k(u)$ given
by
$$
\big(T_k(u)\big)_\mu(x,t):=\mu\int_{-\infty}^t
e^{\mu(s-t)}T_k(u(x,s))\;ds,
$$
extending $T_k(u)$ by $0$ for $s<0$. Observe that $(T_k(u))_\mu\in
L^{p}(0,T;X_0^{s,p}(\Omega))\cap L^\infty(\Omega_T)$, it is differentiable
for a.e. $t\in (0,T)$ with
\begin{eqnarray*}
&&|(T_k(u)\big)_\mu(x,t)|\le k(1-e^{-\mu t})<k\quad \hbox{a.e. in }
\Omega_T,\\[2mm]
&&\frac{\partial(T_k(u))_\mu}{\partial
t}=\mu\big(T_k(u)-(T_k(u))_\mu\big).
\end{eqnarray*}
After computation, we can get
$$
(T_k(u))_\mu\to T_k(u)\quad \hbox{strongly in }
L^{p}(0,T;X_0^{s,p}(\Omega)).
$$

Let us take now a sequence $\{\psi_j\}$ of $C_0^\infty(\Omega)$
functions that strongly converge to $u_0$ in $L^1(\Omega)$, and set
$$
\eta_{\mu,j}(u)\equiv (T_k(u))_{\mu}+e^{-\mu t}T_k(\psi_j).
$$
The definition of $\eta_{\mu,j}$, which is a smooth approximation of
$T_k(u)$, is needed to deal with a nonzero initial datum (see also
\cite{PO}). Note that this function has the following properties:
\begin{eqnarray}
\left\{
\begin{array}{ll}
(\eta_{\mu,j}(u))_t=\mu(T_k(u)-\eta_{\mu,j}(u)),\\[1mm]
\eta_{\mu,j}(u)(0)=T_k(\psi_j),\\[1mm]
|\eta_{\mu,j}(u)|\le k,\\[1mm]
\eta_{\mu,j}(u)\to T_k(u) \quad \textmd{strongly in }
L^{p}(0,T;X_0^{s,p}(\Omega)), \textmd{ as } \mu\to +\infty.
\end{array}
\right.
\end{eqnarray}

Fix a positive number $k$. Let $h>k$. We choose
$$
w_n=T_{2k}\big(u_n-T_h(u_n)+T_k(u_n)-\eta_{\mu,j}(u)\big)
$$
as a test function in (\ref{appro}). Combining the arguments in Step 2 of the proof of Theorem 1.1 in \cite{ZZ} and the Lemma 3.6 in \cite{AAB} together with the nonnegativity and monotonicity of the sequence $\{u_n\}$, we can conclude that
\begin{eqnarray*}
\limsup\limits_{n\to \infty}\int_0^T\!\!\int_{D_\Omega}|T_k(u_n(x,t))-T_k(u_n(y,t))|^p\,d\nu dt\\
\le \int_0^T\!\!\int_{D_\Omega}|T_k(u(x,t))-T_k(u(y,t))|^p\,d\nu dt.
\end{eqnarray*}
It follows from
\begin{eqnarray*}
T_k(u_n)\rightharpoonup T_k(u) \quad\textmd{weakly in } L^p(0,T; X_0^{s,p}(\Omega))
\end{eqnarray*}
that
\begin{eqnarray}
\label{strong convengence}
T_k(u_n)\to T_k(u) \quad\textmd{strongly in } L^p(0,T; X_0^{s,p}(\Omega)).
\end{eqnarray}

\smallskip

\textbf{Step 3.} Show that $u$ is a renormalized solution.

Define the function $G_k(s)=s-T_k(s)$. For given $h>0$, using $T_1(G_h(u_n))$ as a test function in (\ref{appro}), we find
\begin{eqnarray*}
&&\int_{\{|u_n|>h\}} \Theta_1(u_n\mp
h)(T)\,dx-\int_{\{|u_{0n}|>h\}}\Theta_1(u_{0n}\mp h)\,dx \\
&&\quad +\frac12 \int_0^T\!\!\int_{D_\Omega} |U_n(x,y,t)|^{p-2}U_n(x,y,t)\\[2mm]
&&\qquad\qquad\cdot [T_1(G_h(u_n))(x,t)-T_1(G_h(u_n))(y,t)]\,d\nu dt \\
&&\le \int_\Omega f_n T_1(G_h(u_n))\,dxdt,
\end{eqnarray*}
which yields that
\begin{eqnarray*}
&&\frac12 \int_0^T\!\!\int_{D_\Omega} |U_n(x,y,t)|^{p-2}U_n(x,y,t)\\
&&\qquad\quad\cdot [T_1(G_h(u_n))(x,t)-T_1(G_h(u_n))(y,t)]\,d\nu dt\\
&& \le \int_{\{|u_n|>h\}} |f_n|\,dxdt+\int_{\{|u_{0n}|>h\}}|u_{0n}|\,dx.
\end{eqnarray*}
It is not difficult to see that
\begin{eqnarray*}
&&|U_n(x,y,t)|^{p-2}U_n(x,y,t)[T_1(G_h(u_n))(x,t)-T_1(G_h(u_n))(y,t)]\\[2mm]
&&=T'_1(G_h(\xi_n))G'_h(\xi_n)|U_n(x,y,t)|^{p}\ge 0.
\end{eqnarray*}
Recalling the convergence of $\{u_n\}$ in $C([0,T];L^1(\Omega))$, we
have
$$
\lim\limits_{h\to +\infty} {\rm meas}\{(x,t)\in \Omega_T: |u_n|>h\}=0
\quad\textmd{ uniformly with respect to } n.
$$

Since for all $(u_n(x,t), u_n(y,t))\in R_h$,
\begin{eqnarray*}
&&|U_n(x,y,t)|^{p-2}U_n(x,y,t)[T_1(G_h(u_n))(x,t)-T_1(G_h(u_n))(y,t)]\\
&&\ge |U_n(x,y,t)|^{p-1},
\end{eqnarray*}
by using Fatou's lemma and passing to the limit first in $n$ then in $h$, we obtain the renormalized condition
\begin{equation}
\label{renormalized condition}
\lim\limits_{h\to +\infty}\int\!\!\int\!\!\int_{\{(u(x,t), u(y,t))\in R_h\}} |U(x,y,t)|^{p-1}\,d\nu dt=0.
\end{equation}

Let $S\in W^{1, \infty}(\mathbb R)$ be such that ${\rm supp}\,S'\subset [-M,M]$ for some $M>0$. For every $\varphi\in C^1(\bar\Omega_T)$ with $\varphi=0$ in $\mathcal C\Omega\times (0,T)$ and $\varphi(\cdot,T)=0$ in $\Omega$, $S'(u_n)\varphi$ is a test
function in (\ref{appro}). It yields
\begin{eqnarray}
\label{3-18} &&\int_0^T\!\!\int_\Omega \frac{\partial
S(u_n)}{\partial t}\varphi\,dxdt \nonumber \\
&&\quad+\frac 12\int_0^T\!\!\int_{D_\Omega}|U_n(x,y,t)|^{p-2}U_n(x,y,t)(\varphi(x,t)-\varphi(y,t))\nonumber\\
&&\qquad\qquad \cdot \frac{S'(u_n)(x,t)+S'(u_n)(y,t)}{2}\,d\nu dt\nonumber\\
&&\quad+\frac 12\int_0^T\!\!\int_{D_\Omega}|U_n(x,y,t)|^{p-2}U_n(x,y,t)(S'(u_n)(x,t)-S'(u_n)(y,t))\nonumber\\
&&\qquad\qquad\cdot \frac{\varphi(x,t)+\varphi(y,t)}{2}\,d\nu dt\nonumber\\
&&=\int_0^T\!\!\int_\Omega f_nS'(u_n)\varphi\,dxdt.
\end{eqnarray}

First we consider the first term on the left-hand side of
(\ref{3-18}). Since $S$ is bounded and continuous, (\ref{3-4})
implies that $S(u_n)$ converges to $S(u)$ a.e. in $\Omega_T$ and weakly-*
in $L^\infty(\Omega_T)$. Then $\frac{\partial S(u_n)}{\partial t}$
converges to $\frac{\partial S(u)}{\partial t}$ in $D'(\Omega_T)$ as $n\to
+\infty$, that is
$$
\int_0^T\!\!\int_\Omega  S(u_n)\frac{\partial \varphi}{\partial t}\,dxdt\to \int_0^T\!\!\int_\Omega
S(u)\frac{\partial \varphi}{\partial t}\,dxdt.
$$
For the right-hand side of (\ref{3-18}), thanks to the strong
convergence of $f_n$, it is easy to pass to the limit:
$$
\int_0^T\!\!\int_\Omega f_nS'(u_n)\varphi\,dxdt \to \int_0^T\!\!\int_\Omega fS'(u)\varphi\,dxdt, \quad \textmd{as } n\to +\infty.
$$

For the other terms on the left-hand side of (\ref{3-18}), we claim that
\begin{eqnarray*}
&&I_1=\int_0^T\!\!\int_{D_\Omega}|U_n(x,y,t)|^{p-2}U_n(x,y,t)(\varphi(x,t)-\varphi(y,t))\\
&&\qquad\qquad\cdot\frac{S'(u_n)(x,t)+S'(u_n)(y,t)}{2}\,d\nu dt\\
&&\quad \to \int_0^T\!\!\int_{D_\Omega}|U(x,y,t)|^{p-2}U(x,y,t)(\varphi(x,t)-\varphi(y,t))\\
&&\qquad\qquad\cdot \frac{S'(u)(x,t)+S'(u)(y,t)}{2}\,d\nu dt, \quad\textmd{as } n\to +\infty.
\end{eqnarray*}
%

Assume that ${\rm supp}\,S'\subset [-M,M]$. Set
\begin{eqnarray*}
&&D_1=\{(x,y,t)\in D_\Omega\times (0,T): u_n(x,t)\ge M, u_n(y,t)\ge M\},\\
&&D_2=\{(x,y,t)\in D_\Omega\times (0,T): u_n(x,t)\le M, u_n(y,t)\le M\},\\
&&D_3=\{(x,y,t)\in D_\Omega\times (0,T): u_n(x,t)\ge M, u_n(y,t)\le M\},\\
&&D_4=\{(x,y,t)\in D_\Omega\times (0,T): u_n(x,t)\le M, u_n(y,t)\ge M\}.
\end{eqnarray*}
Then $$D_\Omega\times (0,T)=D_1\cup D_2\cup D_3\cup D_4.$$

In $D_1$ we have $$
S'(u_n)(x,t)=S'(u_n)(y,t)=0,
$$ then $I_1=0$.

In $D_2$ we have
$$
u_n(x,t)=T_M(u_n(x,t)), \quad u_n(y,t)=T_M(u_n(y,t)).
$$
From the strong convergence (\ref{strong convengence}), we know that
\begin{eqnarray*}
&&\frac{|T_M(u_n(x,t))-T_M(u_n(y,t))|^{p-2} (T_M(u_n(x,t))-T_M(u_n(y,t)))}{|x-y|^{\frac{(N+ps)(p-1)}{p}}}\\
&&\quad \to \frac{|T_M(u(x,t))-T_M(u(y,t))|^{p-2} (T_M(u(x,t))-T_M(u(y,t)))}{|x-y|^{\frac{(N+ps)(p-1)}{p}}}\\
&&\mbox{strongly in } L^{\frac{p}{p-1}}(D_\Omega\times(0,T)).
\end{eqnarray*}
Moreover, $u_n\to u$ a.e. in $\Omega_T$, $S\in W^{1,\infty}(\mathbb R)$ and $\varphi\in C^1(\overline\Omega_T)$ with $\varphi=0$ in $\mathcal C\Omega\times (0,T)$ imply that
$$
\frac{\varphi(x,t)-\varphi(y,t)}{|x-y|^{\frac{N+ps}{p}}}\in L^p(D_\Omega\times(0,T))
$$
and
\begin{eqnarray*}
&&\frac{\varphi(x,t)-\varphi(y,t)}{|x-y|^{\frac{N+ps}{p}}}\frac{S'(u_n)(x,t)+S'(u_n)(y,t)}{2}\chi_{\{u_n(x,t)\le M,u_n(y,t)\le M\}}\\
&&\quad \rightharpoonup \frac{\varphi(x,t)-\varphi(y,t)}{|x-y|^{\frac{N+ps}{p}}}\frac{S'(u)(x,t)+S'(u)(y,t)}{2}\chi_{\{u(x,t)\le M,u(y,t)\le M\}}\\
&&\mbox{weakly in } L^{p}(D_\Omega\times(0,T)).
\end{eqnarray*}
Thus we have
\begin{eqnarray*}
&&\int\!\!\int_{D_2}|U_n(x,y,t)|^{p-2}U_n(x,y,t)(\varphi(x,t)-\varphi(y,t))\\
&&\qquad\qquad\cdot\frac{S'(u_n)(x,t)+S'(u_n)(y,t)}{2}\,d\nu dt\\
&&\to \int\!\!\int_{\{u(x,t)\le M, u(y,t)\le M\}}|U(x,y,t)|^{p-2}U(x,y,t)(\varphi(x,t)-\varphi(y,t))\\
&&\qquad\qquad\cdot\frac{S'(u)(x,t)+S'(u)(y,t)}{2}\,d\nu dt, \quad\textmd{as } n\to \infty.
\end{eqnarray*}

In $D_3$, if $u_n(x,t)\le M+1$, then it can be done similarly to the estimates in $D_2$. On the other hand, if $u_n(x,t)\ge M+1$, then
$$
\max\{u_n(x,t), u_n(y,t)\}\ge M+1 \quad \mbox{and} \quad  \min\{u_n(x,t),u_n(y,t)\}\le M.
$$
It follows from (\ref{renormalized condition}) that
$$
\lim_{M\to \infty}\lim_{n\to \infty}\int\!\!\int\!\!\int_{\{(u_n(x,t), u_n(y,t))\in R_M\}}|U_n(x,y,t)|^{p-1}\,d\nu dt=0.
$$
Then we observe
\begin{eqnarray*}
&&\!\!\!\!\!\!\!\!\!\!\lim_{M\to \infty}\lim_{n\to \infty}\int\!\!\int_{\{u_n(x,t)\ge M+1, u_n(y,t)\le M\}}|U_n(x,y,t)|^{p-2}U_n(x,y,t)(\varphi(x,t)-\varphi(y,t))\\
&&\quad\qquad\qquad \cdot \frac{S'(u_n)(x,t)+S'(u_n)(y,t)}{2}\,d\nu dt=0.
\end{eqnarray*}
The estimates in $D_4$ can be done similarly.

Therefore, we have
\begin{eqnarray*}
&&\lim_{M\to \infty}\lim_{n\to \infty}I_1\\
&&=\lim_{M\to \infty}\int\!\!\int_{\{u(x,t)\le M, u(y,t)\le M\}}|U(x,y,t)|^{p-2}U(x,y,t)(\varphi(x,t)-\varphi(y,t))\\
&&\qquad\qquad\qquad\cdot \frac{S'(u)(x,t)+S'(u)(y,t)}{2}\,d\nu dt\\
&&\quad+2\lim_{M\to \infty}\int\!\!\int_{\{M\le u(x,t)\le M+1, u(y,t)\le M\}}|U(x,y,t)|^{p-2}U(x,y,t)(\varphi(x,t)-\varphi(y,t))\\
&&\qquad\qquad\qquad\cdot \frac{S'(u)(x,t)+S'(u)(y,t)}{2}\,d\nu dt\\
&&=\int_0^T\!\!\int_{D_\Omega}|U(x,y,t)|^{p-2}U(x,y,t)(\varphi(x,t)-\varphi(y,t))\\
&&\qquad\qquad\qquad\cdot \frac{S'(u)(x,t)+S'(u)(y,t)}{2}\,d\nu dt.
\end{eqnarray*}

The third term on the left-hand side of (\ref{3-18}) can be argued similarly. Therefore, we obtain
\begin{eqnarray*}
&&\int_0^T\!\!\int_\Omega \frac{\partial
S(u)}{\partial t}\varphi\,dxdt \nonumber \\
&&\quad+\frac 12\int_0^T\!\!\int_{D_\Omega}|U(x,y,t)|^{p-2}U(x,y,t)(\varphi(x,t)-\varphi(y,t))\nonumber\\
&&\qquad\qquad\cdot \frac{S'(u)(x,t)+S'(u)(y,t)}{2}\,d\nu dt\nonumber\\
&&\quad+\frac 12\int_0^T\!\!\int_{D_\Omega}|U(x,y,t)|^{p-2}U(x,y,t)(S'(u)(x,t)-S'(u)(y,t))\nonumber\\
&&\qquad\qquad\cdot\frac{\varphi(x,t)+\varphi(y,t)}{2}\,d\nu dt\nonumber\\
&&=\int_0^T\!\!\int_\Omega f S'(u)\varphi\,dxdt,
\end{eqnarray*}
that is
\begin{eqnarray*}
&&-\int_\Omega
S(u_0)\varphi(x,0)\,dx-\int_0^T\!\!\int_\Omega S(u)\frac{\partial
\varphi}{\partial t}\,dxdt \nonumber\\
&&\quad+\frac 12\int_0^T\!\!\int_{D_\Omega} |U(x,y,t)|^{p-2}U(x,y,t)[(S'(u)\varphi)(x,t)-(S'(u)\varphi)(y,t)]\,d\nu dt \nonumber\\
&&=\int_0^T\!\!\int_\Omega fS'(u)\varphi\,dxdt
\end{eqnarray*}
for any $\varphi\in C^1(\bar\Omega_T)$ with $\varphi=0$ in $\mathcal C\Omega\times (0,T)$ and $\varphi(\cdot,T)=0$ in $\Omega$. This completes the proof of the existence of renormalized solutions.

\medskip

(2) Uniqueness of renormalized solutions.

Now we  prove the uniqueness of renormalized solutions for problem
(\ref{main}) by choosing an appropriate test function motivated by
\cite {BMR} and \cite {BW}. Let $u$ and $v$ be two renormalized
solutions for problem (\ref{main}). For
$\sigma>0$, let $S_\sigma$ be the function defined by
\begin{eqnarray}\label{3-19}
\left\{
\begin{array}{ll}
 \displaystyle
S_\sigma(r)=r &\hbox{ if } |r|< \sigma,\\
\displaystyle S_\sigma(r)=(\sigma+\frac 12)\mp\frac
12(r\mp(\sigma+1))^2 &\hbox{ if } \sigma\le
\pm r\le \sigma+1,\\
\displaystyle S_\sigma(r)=\pm(\sigma+\frac 12) &\hbox{ if } \pm
r>\sigma+1.
\end{array}
\right.
\end{eqnarray}
It is obvious that
\begin{eqnarray*}
\left\{
\begin{array}{ll}
\displaystyle
S'_\sigma(r)=1 &\hbox{ if } |r|< \sigma,\\[2mm]
S'_\sigma(r)=\sigma+1-|r| &\hbox{ if } \sigma\le
|r|\le \sigma+1,\\[2mm]
S'_\sigma(r)=0 &\hbox{ if } |r|>\sigma+1.
\end{array}
\right.
\end{eqnarray*}

It is easy to check  $S_\sigma\in W^{1,\infty}(\mathbb R)$ with
${\rm supp}\,S'_\sigma\subset [-\sigma-1,\sigma+1]$.
Therefore,  we may take $S=S_\sigma$ in (\ref{1-5}) to have
\begin{eqnarray*}
&&\int_0^T\!\!\int_\Omega \frac{\partial
S_\sigma(u)}{\partial t}\varphi\,dxdt \nonumber \\
&&\quad+\frac 12\int_0^T \!\!\int_{D_\Omega}|U(x,y,t)|^{p-2}U(x,y,t)(\varphi(x,t)-\varphi(y,t))\nonumber\\
&&\qquad\qquad\cdot\frac{S_\sigma'(u)(x,t)+S_\sigma'(u)(y,t)}{2}\,d\nu dt\nonumber\\
&&\quad+\frac 12\int_0^T \!\!\int_{D_\Omega}|U(x,y,t)|^{p-2}U(x,y,t)(S_\sigma'(u)(x,t)-S_\sigma'(u)(y,t))\nonumber\\
&&\qquad\qquad\cdot\frac{\varphi(x,t)+\varphi(y,t)}{2}\,d\nu dt\nonumber\\
&&=\int_0^T\!\!\int_\Omega f S_\sigma'(u)\varphi\,dxdt
\end{eqnarray*}
and
\begin{eqnarray*}
&&\int_0^T\!\!\int_\Omega \frac{\partial
S_\sigma(v)}{\partial t}\varphi\,dxdt \nonumber \\
&&\quad+\frac 12\int_0^T\!\!\int_{D_\Omega}|V(x,y,t)|^{p-2}V(x,y,t)(\varphi(x,t)-\varphi(y,t))\nonumber\\
&&\qquad\qquad\cdot\frac{S_\sigma'(v)(x,t)+S_\sigma'(v)(y,t)}{2}\,d\nu dt\nonumber\\
&&\quad+\frac 12\int_0^T\!\!\int_{D_\Omega}|V(x,y,t)|^{p-2}V(x,y,t)(S_\sigma'(v)(x,t)-S_\sigma'(v)(y,t))\nonumber\\
&&\qquad\qquad\cdot\frac{\varphi(x,t)+\varphi(y,t)}{2}\,d\nu dt\nonumber\\
&&=\int_0^T\!\!\int_\Omega f S_\sigma'(v)\varphi\,dxdt.
\end{eqnarray*}

For every fixed $k>0$, we plug $\varphi =T_k(S_\sigma(u)-S_\sigma(v))$ as a test function
in the above equalities and subtract them to obtain that
\begin{equation}
\label{3-20} J_0+J_1+J_2=J_3,
\end{equation}
where
\begin{eqnarray*}
&&J_0=\int_0^T \Big\langle\frac{\partial
(S_\sigma(u)-S_\sigma(v))}{\partial t},
T_k(S_\sigma(u)-S_\sigma(v))\Big\rangle\,dt,\\
&&J_1=\frac 12\int_0^T\!\!\int_{D_\Omega}\Big[\frac{S_\sigma'(u)(x,t)+S_\sigma'(u)(y,t)}{2}|U(x,y,t)|^{p-2}U(x,y,t)\\
&&\qquad\qquad\qquad-\frac{S_\sigma'(v)(x,t)+S_\sigma'(v)(y,t)}{2}|V(x,y,t)|^{p-2}V(x,y,t)\Big]\\[2mm]
&&\qquad\qquad\qquad\cdot [T_k(S_\sigma(u)-S_\sigma(v))(x,t)-T_k(S_\sigma(u)-S_\sigma(v))(y,t)]\,d\nu dt,\\[2mm]
&&J_2=\frac 12\int_0^T\!\!\int_{D_\Omega}\big[|U(x,y,t)|^{p-2}U(x,y,t)\cdot (S_\sigma'(u)(x,t)-S_\sigma'(u)(y,t))\\
&&\qquad\qquad\qquad-|V(x,y,t)|^{p-2}V(x,y,t)\cdot (S_\sigma'(v)(x,t)-S_\sigma'(v)(y,t))\big]\\[2mm]
&&\qquad\qquad\qquad\cdot \frac{T_k(S_\sigma(u)-S_\sigma(v))(x,t)+T_k(S_\sigma(u)-S_\sigma(v))(y,t)}{2}\,d\nu dt,\\[2mm]
&&J_3=\int_0^T\!\!\int_\Omega
f(S'_\sigma(u)-S'_\sigma(v))T_k(S_\sigma(u)-S_\sigma(v))\,dxdt.
\end{eqnarray*}

We estimate $J_0$, $J_1$, $J_2$ and $J_3$ one by one. Recalling the
definition of $\Theta_k(r)$, $J_0$ can be written as
\begin{eqnarray*}
J_0=\int_\Omega \Theta_k(S_\sigma(u)-S_\sigma(v))(T)\,dx-\int_\Omega
\Theta_k(S_\sigma(u)-S_\sigma(v))(0)\,dx.
\end{eqnarray*}
Due to the same initial condition for $u$ and $v$, and the
properties of $\Theta_k$, we get
$$
J_0=\int_\Omega \Theta_k(S_\sigma(u)-S_\sigma(v))(T)\,dx\ge 0.
$$

Writing
\begin{eqnarray*}
&&J_1=\frac 12\int_0^T\!\!\int_{D_\Omega} (|U(x,y,t)|^{p-2}U(x,y,t)-|V(x,y,t)|^{p-2}V(x,y,t))\\[2mm]
&&\qquad\qquad\cdot[T_k(S_\sigma(u)-S_\sigma(v))(x,t)-T_k(S_\sigma(u)-S_\sigma(v))(y,t)]\,d\nu dt\\[2mm]
&&\qquad+\frac 12\int_0^T\!\!\int_{D_\Omega} \left(1-\frac{S_\sigma'(u)(x,t)+S_\sigma'(u)(y,t)}{2}\right)|U(x,y,t)|^{p-2}U(x,y,t)\\[2mm]
&&\qquad\qquad\cdot[T_k(S_\sigma(u)-S_\sigma(v))(x,t)-T_k(S_\sigma(u)-S_\sigma(v))(y,t)]\,d\nu dt\\[2mm]
&&\qquad+\frac 12\int_0^T\!\!\int_{D_\Omega} \left(\frac{S_\sigma'(v)(x,t)+S_\sigma'(v)(y,t)}{2}-1\right)|V(x,y,t)|^{p-2}V(x,y,t)\\[2mm]
&&\qquad\qquad \cdot[T_k(S_\sigma(u)-S_\sigma(v))(x,t)-T_k(S_\sigma(u)-S_\sigma(v))(y,t)]\,d\nu dt\\[2mm]
&&\quad:=J^1_1+J^2_1+J^3_1,
\end{eqnarray*}
and setting  $\sigma\ge k$, we have
\begin{eqnarray}\label{J11}
&&\!\!\!\!\!\!\!\!\!\!J^1_1\ge \frac 12\int\!\!\int\!\!\int_{\{|u-v|\le k\}\cap\{|u|,|v|\le k\}}(|U(x,y,t)|^{p-2}U(x,y,t)-|V(x,y,t)|^{p-2}V(x,y,t)) \nonumber\\[2mm]
&&\qquad\qquad \cdot[U(x,y,t)-V(x,y,t)]\,d\nu dt.
\end{eqnarray}
By the Lebesgue dominated convergence theorem, we conclude that
$$
J_1^2, J_1^3\to 0, \quad\textmd{as } \sigma\to +\infty.
$$
Furthermore, we have
\begin{eqnarray*}
|J_2|&\le& C\Big(\int\!\!\int\!\!\int_{\{(u(x,t), u(y,t))\in R_\sigma\}} |U(x,y,t)|^{p-1}\,d\nu dt\\
&&+\int\!\!\int\!\!\int_{\{(v(x,t), v(y,t))\in R_\sigma\}}|V(x,y,t)|^{p-1}\,d\nu dt\Big).
\end{eqnarray*}
From the above estimates and (i) in Definition 1.1, we obtain
$$
\lim\limits_{\sigma\to +\infty} (|J^2_1|+|J^3_1|+|J_2|)=0.
$$

Observing
\begin{eqnarray*}
f(S'_\sigma(u)-S'_\sigma(v))\to 0 \quad\mbox{strongly in } L^1(\Omega_T)
\end{eqnarray*}
as $\sigma\to +\infty$ and using the Lebesgue dominated convergence
theorem, we deduce that
$$
\lim\limits_{\sigma\to +\infty}|J_3|=0.
$$

Therefore, sending $\sigma\to +\infty$ in (\ref{3-20}) and recalling
(\ref{J11}), we have
\begin{eqnarray*}
\int\!\!\int\!\!\int_{\{|u|\le \frac k2,|v|\le \frac k2\}}(|U(x,y,t)|^{p-2}U(x,y,t)-|V(x,y,t)|^{p-2}V(x,y,t)) \nonumber\\[2mm]
\cdot[U(x,y,t)-V(x,y,t)]\,d\nu dt=0,
\end{eqnarray*}
which implies $U=V$ a.e. on the set $\big\{|u|\le
\frac k2,|v|\le \frac k2\big\}$. Since $k$ is arbitrary, we conclude
that
$$
U=V \quad\textmd{for a.e. } x,y\in \mathbb R^N, t\in [0,T].
$$
It follows from the Poincar\'{e} inequality with $p=1$ in (\ref{poincare}) that
\begin{eqnarray*}
&&\int_0^T\!\!\int_\Omega |u(x,t)-v(x,t)|\,dxdt\\
&&\le C\int_0^T\!\!\int_{D_\Omega} |U(x,y,t)-V(x,y,t)|\,d\nu_1 dt=0,
\end{eqnarray*}
where
$$
d\nu_1=\frac{dxdy}{|x-y|^{N+s}}.
$$
Thus we have $u=v$ a.e. in $\Omega_T$. This completes the proof of Theorem 1.1.
$\quad\Box$

\medskip

Next, we prove that the renormalized solution $u$ is also an entropy
solution of problem (\ref{main}) and the entropy solution of problem
(\ref{main}) is unique.

\medskip

\noindent{\bf Proof of Theorem 1.2.} (1) The renormalized solution is an entropy solution.

Now we choose $v_n=T_k(u_n-\phi)$ as a test function in
(\ref{appro}) for $k>0$ and  $\phi\in C^1(\bar \Omega_T)$ with
$\phi=0$ in $\mathcal C\Omega\times (0,T)$. Following the arguments in \cite{AABP}, we can prove the existence of entropy
solutions.

\medskip (2) Uniqueness of entropy solutions.

Suppose that $u$ and $v$ are two entropy solutions of problem
(\ref{main}). Let $\{u_n\}$ be a sequence constructed in
(\ref{appro}). Choosing $S_\sigma(u_n)$ as a test function in (\ref{1-4}) for
entropy solution $v$, we have
\begin{eqnarray}
\label{esti4} &&\int_\Omega
\Theta_k(v-S_\sigma(u_n))(T)\,dx-\int_\Omega
\Theta_k(u_0-S_\sigma(u_{0n}))\,dx \nonumber\\
&&\quad+\int_0^T \langle (u_n)_t, S_\sigma'(u_n)T_k(v-S_\sigma(u_n)
\rangle\,dt\nonumber
\\
&&\quad+\frac 12\int_0^T\!\!\int_{D_\Omega} |V(x,y,t)|^{p-2}V(x,y,t)[T_k(v(x,t)-S_\sigma(u_n)(x,t)) \nonumber\\
&&\qquad\qquad-T_k(v(y,t)-S_\sigma(u_n)(y,t))]\,d\nu dt \nonumber \\
&&\le \int_0^T\!\!\int_\Omega f T_k(v-S_\sigma(u_n))\,dxdt.
\end{eqnarray}

In order to deal with the third term on the left-hand side of
(\ref{esti4}), we take $S_\sigma'(u_n)\Psi$ with
$\Psi=T_k(v-S_\sigma(u_n))$ as a test function for problem
(\ref{appro}) to obtain
\begin{eqnarray}
\label{esti5} &&\int_0^T \langle(u_n)_t, S_\sigma'(u_n)\Psi
\rangle\,dt\nonumber\\
&&\quad+\frac 12\int_0^T\!\!\int_{D_\Omega}|U_n(x,y,t)|^{p-2}U_n(x,y,t)(S_\sigma'(u_n)(x,t)-S_\sigma'(u_n)(y,t))\nonumber\\
 &&\qquad\qquad\cdot \frac{\Psi(x,t)+\Psi(y,t)}{2}\,d\nu dt\nonumber\\
&&\quad+\frac 12\int_0^T\!\!\int_{D_\Omega}|U_n(x,y,t)|^{p-2}U_n(x,y,t)(\Psi(x,t)-\Psi(y,t))\nonumber\\
 &&\qquad\qquad\cdot \frac{S_\sigma'(u_n)(x,t)+S_\sigma'(u_n)(y,t)}{2}\,d\nu dt\nonumber\\
&&=\int_0^T\!\!\int_\Omega f_n S_\sigma'(u_n)\Psi\,dxdt.
\end{eqnarray}
Thus we deduce from (\ref{esti4}) and (\ref{esti5}) that
\begin{eqnarray*}
&&\int_\Omega \Theta_k(v-S_\sigma(u_n))(T)\,dx-\int_\Omega
\Theta_k(u_0-S_\sigma(u_{0n}))\,dx\\
&&\quad-\frac 12\int_0^T \!\!\int_{D_\Omega}|U_n(x,y,t)|^{p-2}U_n(x,y,t)(S_\sigma'(u_n)(x,t)-S_\sigma'(u_n)(y,t))\nonumber\\
 &&\qquad\qquad\cdot \frac{\Psi(x,t)+\Psi(y,t)}{2}\,d\nu dt\\
&&\quad-\frac 12\int_0^T \!\!\int_{D_\Omega}|U_n(x,y,t)|^{p-2}U_n(x,y,t)(\Psi(x,t)-\Psi(y,t))\nonumber\\
 &&\qquad\qquad\cdot \frac{S_\sigma'(u_n)(x,t)+S_\sigma'(u_n)(y,t)}{2}\,d\nu dt\\
&&\quad+\frac 12\int_0^T\!\!\int_{D_\Omega} |V(x,y,t)|^{p-2}V(x,y,t)(\Psi(x,t)-\Psi(y,t))\,d\nu dt\\
&&\le \int_0^T\!\!\int_\Omega  f
T_k(v-S_\sigma(u_n))\,dxdt-\int_0^T\!\!\int_\Omega f_n
S_\sigma'(u_n)T_k(v-S_\sigma(u_n))\,dxdt.
\end{eqnarray*}

We will pass to the limit as $n\to +\infty$ and $\sigma\to +\infty$
successively. Let us denote $A_3$ for the third term on the
left-hand side of the above equality for simplicity. Recalling the definition of $S_\sigma'$,
we have
$$
|A_3|\le k\int\!\!\int\!\!\int_{\{(u_n(x,t), u_n(y,t))\in R_\sigma\}}|U_n(x,y,t)|^{p-1}\,d\nu dt.
$$

Observe that
\begin{eqnarray*}
&&\int_0^T\!\!\int_{D_\Omega} |V(x,y,t)|^{p-2}V(x,y,t)(\Psi(x,t)-\Psi(y,t))\,d\nu dt\\
&&\quad-\int_0^T \!\!\int_{D_\Omega}|U_n(x,y,t)|^{p-2}U_n(x,y,t)(\Psi(x,t)-\Psi(y,t)) \\
 &&\qquad\qquad\cdot\frac{S_\sigma'(u_n)(x,t)+S_\sigma'(u_n)(y,t)}{2}\,d\nu dt\\
&&=\int_0^T\!\!\int_{D_\Omega} (|V(x,y,t)|^{p-2}V(x,y,t)-|U_n(x,y,t)|^{p-2}U_n(x,y,t))\\
 &&\qquad\qquad\cdot(\Psi(x,t)-\Psi(y,t))\,d\nu dt\\
&&\quad+\int_0^T\!\!\int_{D_\Omega}|U_n(x,y,t)|^{p-2}U_n(x,y,t)(\Psi(x,t)-\Psi(y,t)) \\
 &&\qquad\qquad\cdot\left(1-\frac{S_\sigma'(u_n)(x,t)+S_\sigma'(u_n)(y,t)}{2}\right)\,d\nu dt.
\end{eqnarray*}

Using the similar arguments as in Theorem 1.1 and the Lebesgue dominated convergence
theorem, letting $n\to +\infty$, we obtain
\begin{eqnarray}
\label{3-30} &&\int_\Omega
\Theta_k(v-S_\sigma(u))(T)\,dx-\int_\Omega
\Theta_k(u_0-S_\sigma(u_{0}))\,dx \nonumber\\
&&\quad+\frac 12\int_0^T\!\!\int_{D_\Omega} (|V(x,y,t)|^{p-2}V(x,y,t)-|U(x,y,t)|^{p-2}U(x,y,t)\\
 &&\qquad\qquad\cdot[T_k(v-S_\sigma(u))(x,t)-T_k(v-S_\sigma(u))(y,t)]\,d\nu dt
\nonumber\\
&&\le \int_0^T\!\!\int_\Omega  f(1-S_\sigma'(u))
T_k(v-S_\sigma(u))\,dxdt \nonumber \\
&&\quad+k\int\!\!\int\!\!\int_{\{(u(x,t), u(y,t))\in R_\sigma\}}|U(x,y,t)|^{p-1}\,d\nu dt\nonumber\\
&&\quad +\frac 12 \int_0^T \!\!\int_{D_\Omega}|U(x,y,t)|^{p-2}U(x,y,t)[T_k(v-S_\sigma(u))(x,t)-T_k(v-S_\sigma(u))(y,t)] \nonumber\\
&&\qquad\qquad\cdot\left(\frac{S_\sigma'(u)(x,t)+S_\sigma'(u)(y,t)}{2}-1\right)\,d\nu dt.
\end{eqnarray}


Now we let $\sigma\to +\infty$. Since
$$
|\Theta_k(v-S_\sigma(u))(T)|\le k(|v(T)|+|u(T)|),\quad
|\Theta_k(u_0-S_\sigma(u_0))|\le k|u_0|,
$$
by the Lebesgue dominated convergence theorem, we have
$$
\int_\Omega \Theta_k(u_0-S_\sigma(u_{0}))\,dx\to 0, \quad \int_\Omega
\Theta_k(v-S_\sigma(u))(T)\,dx\to \int_\Omega \Theta_k(v-u)(T)\,dx.
$$
According to the fact that
$$
\lim\limits_{\sigma\to +\infty}\int\!\!\int\!\!\int_{\{(u(x,t), u(y,t))\in R_\sigma\}}|U(x,y,t)|^{p-1}\,d\nu dt=0
$$
and Fatou's lemma, we deduce from (\ref{3-30}) that
\begin{eqnarray*}
&&\int_\Omega \Theta_k(v-u)(T)\,dx \\
&&\quad+\frac 12\int_0^T\!\!\int\!\!\int_{\{|u|\le \frac k2,|v|\le \frac k2\}} (|V(x,y,t)|^{p-2}V(x,y,t)-|U(x,y,t)|^{p-2}U(x,y,t)\\
 &&\qquad\qquad\cdot[V(x,y,t)-U(x,y,t)]\,d\nu dt\le 0.
\end{eqnarray*}
Using the positivity of $\Theta_k$,
we conclude that $u=v$ a.e. in $\Omega_T$. Therefore we obtain the
uniqueness of entropy solutions. This completes the proof of Theorem
1.2. $\quad\Box$

\bigskip
\noindent{\bf Proof of Theorem 1.3.}
First, we suppose that $u_0, v_0\in L^2(\Omega)$ and $f, g\in
L^{p'}(0,T; X_0^{s,p}(\Omega)^*)$. Then by an approximation argument, we can obtain two weak solutions $u$ and $v$  for problems (\ref{main}) and
\begin{eqnarray}
\left\{
\begin{array}{ll}
\displaystyle v_t+(-\Delta)_p^s v=f
&\textmd{in} \quad \Omega_T,\\[2mm]
v=0 &\textmd{in} \quad \mathcal C\Omega\times (0,T),\\[2mm]
v(x,0)=v_0(x) &\textmd{in} \quad \Omega.
\end{array}
\right. \label{v}
\end{eqnarray}
Making use of the approximation argument, we choose
$(u-v)^+\chi_{(0,t)}$ as a test function and subtract the resulting
equalities to get
\begin{eqnarray*}
&&\int_0^t\!\!\int_\Omega (u-v)_t
(u-v)^+\,dxd\tau+\int_0^t\langle(-\Delta)_p^s u-(-\Delta)_p^s v, (u-v)^+\rangle\,d\tau\\
&&=\int_0^t\!\!\int_\Omega(f-g)(u-v)^+\,dxd\tau\le 0.
\end{eqnarray*}
Moreover, from the nonnegativity of the second term in the equality above, we have
\begin{eqnarray*}
&&\frac 12 \int_0^t\!\!\int_\Omega \frac
{d}{dt}[(u-v)^+]^2\,dxd\tau\\
&&=\frac 12 \int_\Omega
[(u-v)^+]^2(t)\,dx-\frac 12 \int_\Omega [(u_0-v_0)^+]^2\,dx\le 0.
\end{eqnarray*}
Recalling $u_0\le v_0$, we conclude that
$$
(u-v)^+=0 \quad\textmd{a.e. in } \Omega_T.
$$
Thus we obtain $u\le v$ a.e. in $\Omega_T$.

Now we consider $u$ and $v$ as the entropy solution (renormalized solution) of problems
(\ref{main}) and (\ref{v}) with $L^1$ data. Find four sequences of
functions $\{f_n\}, \{g_n\}\subset C^\infty_0(\Omega_T)$ and
$\{u_{0n}\}, \{v_{0n}\}\subset C^\infty_0(\Omega)$ strongly
converging respectively to $f, g$ in $L^1(\Omega_T)$ and to $u_0,
v_0$ in $L^1(\Omega)$ such that
\begin{eqnarray*}
&f_n\le g_n, & u_{0n}\le v_{0n},\\
&\|f_n\|_{L^1(\Omega_T)}\le \|f\|_{L^1(\Omega_T)},&
\|g_n\|_{L^1(\Omega_T)}\le \|g\|_{L^1(\Omega_T)},\\
&\|u_{0n}\|_{L^1(\Omega)}\le \|u_0\|_{L^1(\Omega)}, &
\|v_{0n}\|_{L^1(\Omega)}\le \|v_0\|_{L^1(\Omega)}.
\end{eqnarray*}
Thus we use Theorem 1.1 (Theorem 1.2) to construct two approximation
sequences $\{u_n\}$ and $\{v_n\}$ of entropy solutions (renormalized solutions) $u$ and $v$,
and apply the comparison result  above to obtain  $u_n\le v_n$ a.e.
in $\Omega_T$. Moreover, by the uniqueness of entropy solutions (renormalized solutions), we
know $u_n\to u$ and $v_n\to v$ a.e. in $\Omega_T$. Therefore, we
conclude that $u\le v$ a.e. in $\Omega_T$. This completes the proof
of Theorem 1.3. $\quad\Box$

\section{Extensions}
\label{} \setcounter{equation}{0}

In order to fix the ideas and to avoid unessential technicalities, we limited ourselves to the
equations of principal type as the one considered in (\ref{main}). Indeed, inspired by \cite{BPR}, the existence and uniqueness result of nonnegative renormalized solutions obtained in Theorem 1.1 still holds for the following more general nonlinear parabolic equations
\begin{eqnarray}
\left\{
\begin{array}{ll}
\displaystyle \frac{\partial b(u)}{\partial t}-\mathcal L_p u=f
&\textmd{in} \quad \Omega_T,\\[3mm]
u=0 &\textmd{in} \quad \mathcal C\Omega \times (0,T),\\[2mm]
b(u)(x,0)=b(u_0)(x) &\textmd{in} \quad \Omega,
\end{array}
\right. \label{main2}
\end{eqnarray}
where $u_0$ is a nonnegative measurable function such that $b(u_0)\in L^1(\Omega)$, $0\le f\in L^1(\Omega_T)$, $-\mathcal L_p$ is a non-local operator defined by
\begin{eqnarray*}
-\mathcal L_p u(x,t):={\rm P.V.} \int_{\mathbb R^N} |u(x,t)-u(y,t)|^{p-2}(u(x,t)-u(y,t))K(x,y)\,dy,
\end{eqnarray*}
where $(x, t)\in \mathbb R^N\times \mathbb R^+$, and $b: \mathbb R\to \mathbb R$ is a strictly increasing $C^1$-function satisfying that
$$
0<b_0\le b'(s)\le b_1, \quad b(0)=0.
$$
Finally, the kernel $K: \mathbb R^N\times \mathbb R^N\to \mathbb R$ is assumed to be measurable, and satisfies
the following ellipticity/coercivity properties:
$$
\frac{1}{\Lambda|x-y|^{N+sp}}\le K(x,y)\le \frac{\Lambda}{|x-y|^{N+sp}}, \quad \forall x,y\in \mathbb R^N, x\not=y, \Lambda\ge 1,
$$
where $0<s<1<p<N$ such that $ps<N$.

\medskip

The definition of renormalized solutions for problem (\ref{main2}) is as follows.

\smallskip

\noindent{\bf Definition 4.1.}
 \label{def4-1}
A function $u$ defined on $\mathbb R^N\times (0,T]$ is a renormalized solution to problem
{\rm(\ref{main2})} if $b(u)\in
C([0,T];L^1(\Omega))$, $T_k(b(u))\in L^p(0,T;X_0^{s,p}(\Omega))$ for any $k\ge 0$, and the following conditions are satisfied:

\smallskip

(i) $$\displaystyle \lim\limits_{h\to \infty}\int\!\!\int\!\!\int_{\{(x,y,t):(b(u)(x,t), b(u)(y,t))\in R_h\}}|U(x,y,t)|^{p-1}\,d\nu dt=0,$$ where
$$R_h=\Big\{(u,v)\in \mathbb R^{2}: h+1\le \max\{|u|, |v|\} \textmd{ and } (\min\{|u|, |v|\}\le h \textmd{ or } uv<0)\Big\}.$$

(ii) For every function $\varphi\in C^1(\bar\Omega_T)$ with $\varphi=0$ in $\mathcal C\Omega\times (0,T)$ and $\varphi(\cdot,T)=0$ in $\Omega$, and $S\in W^{1,\infty}(\mathbb R)$ which is
piecewise $C^1$ satisfying that $S'$ has a compact support,
\begin{eqnarray*}
&&\!\!\!\!\!\!\!\!\!\!-\int_\Omega
S(b(u_0))\varphi(x,0)\,dx-\int_0^T\!\!\int_\Omega S(b(u))\frac{\partial
\varphi}{\partial t}\,dxdt \nonumber\\
&&\!\!\!\!\!\!\!\!\!+\frac 12\int_0^T\!\!\int_{D_\Omega}  |\tilde U(x,y,t)|^{p-2}\tilde U(x,y,t)[(S'(b(u))\varphi)(x,t)-(S'(b(u))\varphi)(y,t)]\,d\nu dt \nonumber\\
&&\!\!\!\!\!=\int_0^T\!\!\int_\Omega fS'(b(u))\varphi\,dxdt
\end{eqnarray*}
holds, where
$$
\tilde U(x,y,t)=b(u)(x,t)-b(u)(y,t).
$$

To the best of our knowledge, it is an open problem to show the well-posedness of entropy solutions and the equivalence between renormalized and entropy solutions to the general problem (\ref{main2}).

\section*{Acknowledgements}
The authors wish to thank Dr. Xia Zhang for careful reading
an early version of this paper, pointing out a mistake in the proof of Theorem 1.1 and helping them correct the mistake.

K. Teng was supported by the NSFC (No. 11501403) and the Shanxi Province Science Foundation for Youths (No. 2013021001-3). C. Zhang was supported by the NSFC (No. 11671111) and Heilongjiang Province Postdoctoral Startup Foundation (LBH-Q16082). S. Zhou was supported by the NSFC (No. 11571020).





\bibliographystyle{elsarticle-num}
\bibliography{<your-bib-database>}



\end{document}